\documentclass[]{article}
\usepackage{amsmath,amssymb,amsfonts,amsthm,bbm}
\numberwithin{equation}{section}
\numberwithin{table}{section}
\numberwithin{figure}{section}
\usepackage{cite}
\usepackage[linesnumbered,ruled,algosection]{algorithm2e}
\usepackage{caption}
\usepackage{float}
\usepackage{graphicx}
\usepackage{epstopdf}
\usepackage{color,soul}
\def\rr{{\mathbb R}}
\def\EE{\mathbb E}
\def\nn{{\mathbb N}}
\providecommand{\keywords}[1]{\textbf{\textit{\small Keywords }} #1}

\title{Commodity Resource Valuation And Extraction: A Pathwise Programming Approach}
\author{
  Juri Hinz\\
  \and
  Tanya Tarnopolskaya\\
  \and
  Jeremy Yee* \\
}
\date{\today}

\begin{document}
\maketitle
\begin{abstract}
  Complexity and uncertainty associated with commodity resource
  valuation and extraction requires stochastic control methods suitable
  for high dimensional states. Recent progress in duality and
  trajectory-wise techniques has introduced a variety of fresh ideas to
  this field with surprising results. This paper presents a first
  application of this promising development to commodity extraction
  problems. We introduce efficient algorithms for obtaining approximate
  solutions along with a diagnostic technique, which provides a
  quantitative measure for solution performance in terms of the
  distance between the approximate and the optimal control policy.
\end{abstract}

\section{Introduction \label{sec-intro}}
Extraction projects for commodities and their valuation can often be
formulated as optimal stochastic control problems of the switching
type. Such problems usually have no known closed-form solutions and
hence numerical methods often remain the only practical approach to
address important operational and investment questions.  This paper
presents a family of novel methods for calculating approximate
numerical solutions based on primal and dual techniques. While the
primal methods are based on the so-called {\it convex stochastic
  switching} (CSS) and deliver an approximate solution, the dual
methods utilize recent advances of {\it pathwise dynamic programming}
to provide solution improvement and its quality assessment.  This
paper will demonstrate our algorithms by considering the optimal
decision policy and value of natural resource assets associated with
commodity extraction. Numerical results were obtained using the
authors' software package implemented in the statistical language
\emph{R}.  Similar optimal switching problems were considered by
\cite{bayraktar_egami2010, boomsma_etal2012, brennan_schwartz1985,
  carmona_ludkovski2010, tsekrekos_etal2012,
  tsekrekos_yannacopoulos2016}.  This paper will depart from the
standard literature and consider the presence of \emph{mean reversion}
in the commodity price, operational effiencies in the form of
\emph{wasteful extraction}, \emph{physical delivery} obligations to
clients and commodity prices subject to \emph{stochastic volatility}.
These considerations add realism to our study and sheds practical
insights.  

In the next section, we introduce the problem setting. Section
\ref{sec-algo} outlines our numerical approach while Section
\ref{sec-result} contains the numerical results. Section
\ref{sec-end} concludes.

\section{Commodity Extraction As A Stochastic 
  Switching Problem \label{sec-problem}}

This paper will examine the management of a commodity resource with a finite
amount of commodity as studied by \cite{brennan_schwartz1985}.  In this
setting, the decision maker aims to maximize the total expected profit
using controls whose costs are given in Table \ref{bs_table}.  Doing
so, the controller dynamically switches the operational mode of the
resource between \{Opened, Closed, Abandoned\}. When opened, commodity is
extracted, sold at the spot market and a revenue is realized based on
the commodity price.  While closed, commodity is preserved but a maintenance
cost is incurred. The abandoned mode incurs no cost and provides no
revenue. Switching to the abandoned mode causes the mode to remain
permanently there. There is a switching cost between modes.

\begin{table} [h]
  \centering
  \caption{Hypothetical Commodity Resource. 
    \label{bs_table}}
  \begin{tabular}{ll}
    \hline
    Output rate: 10 million pounds/year &
    Inventory level: 150 million pounds \\
    Costs of production: \$0.50/pound &
    Opening/closing costs: \$200,000 \\
    Maintenance costs: \$500,000/year &
    Inflation rate: 8\%/year \\
    Convenience yield: 1\%/year &
    Price variance: 8\%/year \\
    Real estate tax: 2\%/year  &
    Income tax: 50\% \\
    Royalty tax: 0\% &
    Interest rate: l0\%/year \\
    \hline
  \end{tabular}
\end{table}

The original problem setting in \cite{brennan_schwartz1985} was
considered in a continuous time setting and solved using numerical
partial differential equations. However, this paper will consider the
problem as a \emph{Markov decision process} (see \cite{puterman})) in
discrete time and a finite time horizon. The rates, costs and taxes
above will be adjusted accordingly.

Given a finite time horizon $\{0, 1, \dots, T\} \subset \nn$, consider
a controlled Markovian process $(X_t)_{t=0}^T := (P_t, Z_t)_{t=0}^T$
which consists of two parts. The discrete component $(P_t)_{t=0}^T$
describes the evolution of a finite-state controlled Markov chain
which takes values in a finite set $\mathbf{P}$.  Further assume that
at any time $t=0, \dots, T-1$ the controller takes an action $a \in
\mathbf{A}$ from a finite set $\mathbf{A}$ of all admissible actions
in order to cause the one-step transition from the mode $p \in
\mathbf{P}$ to the mode $p' \in \mathbf{P}$ with a probability
$\alpha_{p, p'}(a)$, where $(\alpha_{p, p'}(a))_{p, p' \in
\mathbf{P}}$ are pre-specified stochastic matrices for all $a \in
\mathbf{A}$.  Let us now turn to the evolution of the other component
$(Z_t)_{t=0}^T$ of the state process $(X_t)_{t=0}^T$.  Here, we assume
that it follows an uncontrolled evolution in the Euclidean space
$\rr^d$ driven as
$$
Z_{t+1}=W_{t+1} Z_{t}, \qquad t=0, \dots, T-1
$$
by independent {\it disturbance matrices} $(W_{t})_{t=1}^T$. That is,
the transition kernels ${\cal K}^a_t$ governing the evolution of our
controlled Markov process $(X_t)_{t=0}^T := (P_t, Z_t)_{t=0}^T$ from
time $t$ to $t+1$ is given for each $a \in \mathbf{A}$ by
$$
{\cal K}_t v(p,z)= \sum_{p' \in \mathbf{P}} \alpha_{p, p'}(a)
\EE(v(p', W_{t+1}z)), \quad p \in \mathbf{P}, \, z \in \rr^d, \, t=0,
\dots, T-1
$$
which acts on each function $v: \mathbf{P} \times \rr^d \to \rr$ where
the above expectations are well-defined.

In this work, we use the discrete component $(P_t)_{t=0}^T$ to convey
information regarding the remaining level of commodity and the current
operational mode. More precisely, set $ \mathbf{P}=\{0, 1,
\dots, I\} \times \{1, 2\} $ to capture all possible positions with
the interpretation that the first component $p^{(1)}$ of the position
$(p^{(1)}, p^{(2)}) \in \mathbf{P}$ describes the remaining commodity level
such that $p^{(1)}=0$ represent an exhausted or abandoned resource. The
second component $p^{(2)}$ indicates whether the current operational
mode is closed ($p^{(2)}=1$) or opened ($p^{(2)}=2$).  To introduce
the mode control, we set $\mathbf{A} = \{0, 1, 2\}$ where $a
= 0, 1, 2$ stands for the action to abandon, to close, and to open the
asset, respectively. The second component $(Z_t)_{t=0}^T$ will be used
to capture the market dynamics of the the commodity price. The exact form
of the disturbance matrices depends on the particular choice of the price
process.  If the system
is in the state $(p, z)$, the costs of applying action $a \in
\mathbf{A}$ at time $t=0, \dots, T-1$ are expressed through $r_t(p, z,
a)$. Having arrived at time $t=T$ in the state $(p, z)$, a final {\it
scrap value} $r_T(p, z)$ is collected. Thereby the reward and scrap
functions
$$ r_t: \mathbf{P} \times \rr^d \times \mathbf{A} \to \rr, \quad 
 r_T: \mathbf{P} \times \rr^d \to \rr
$$
are exogenously given for $t=0, \dots, T-1$.  In our applications, the
scrap value is $$r_T(p,z)=\max_{a \in \mathbf{A}}r_T(p, z, a), \quad
(p, z) \in \mathbf{P} \times \rr^d$$ where the reward functions are
specified in the following manner. If the asset is abandoned, there
are neither costs nor revenue
 \[ r_{t}((0, p^{(2)}, z), a)=0, \qquad a \in \mathbf{A}, z \in
\mathbb{R}^d. \] For the case where $p^{(1)}>0$, we define
 \[ r_{t}((p^{(1)}, p^{(2)}, z), a)=h_{t}(z)1_{\{2\}} (a) + m_{t} 1_{
\{1\}} (a) + c_t 1_{\{1,2\}}(a)|p^{(2)}-a| \] with the following
interpretation: If the resource is opened, a revenue is collected which
depends on continuous state component $z$ through a pre-specified
convex function $h_{t}$. If the resource is closed, then non-stochastic
maintenance fee $m_{t}$ is to be paid. Finally, a deterministic
switching fee $c_{t}$ is incurred whenever the operational mode
transitions from opened to closed or vice versa.  At each time $t=0,
\dots, T$ the \emph{decision rule} $\pi_{t}$ is given by a mapping
$\pi_{t}: \mathbf{P} \times \rr^d \to \mathbf{A}$, prescribing at time
$t$ an action $\pi_{t}(p,z) \in \mathbf{A}$ for a given state $(p, z)
\in \mathbf{P} \times \rr^d$.  A sequence $\pi =
(\pi_{t})_{t=0}^{T-1}$ of decision rules is called a \emph{policy}.
For each policy $\pi = (\pi_{t})_{t=0}^{T-1}$, the so-called policy
value $v^{\pi}_{0}(p_0, z_{0})$ is defined as the total expected
reward
\[ v^{\pi}_{0}(p_0, z_{0})=\mathbb{E}^{ (p_0, z_0), \pi } \left[
\sum_{t=0}^{T} r_{t}(P_t, Z_{t}, \pi_{t}(P_t, Z_t)) + r_{t}(P_t,
Z_{t}) \right].
\] 
In this formula $\mathbb{E}^{ (p_0, z_0), \pi }$ stands for the
expectation with respect to the probability distribution of
$(X_t)_{t=0}^T := (P_t, Z_t)_{t=0}^T$ induced by the Markov
transitions from $(P_{t}, Z_{t})$ to $(P_{t+1}, Z_{t+1})$ induced by
the kernels ${\cal K}^{\pi_t(P_t, Z_t)}_t$ for $t=0, \dots, T-1$,
given the initial value $(P_0, Z_0)=(p_0, z_0)$.

Now we turn to the optimization goal.  A policy
$\pi^{*}=(\pi^{*}_{t})_{t=0}^{T-1}$ is called optimal if it maximizes
the total expected reward over all policies $\pi \mapsto
v^{\pi}_{0}(p, z)$. To obtain such policy, one introduces for $t = 0,
\dots, T-1$ the so-called {\it Bellman operator}
\begin{equation} 
{\cal T}_{t}v(p,z)=\max_{a \in \mathbf{A}} \left[
r_{t}(p,z , a) + \sum_{p'\in\mathbf{P}}
\alpha_{p,p'}(a)\mathbb{E}[v(p',W_{t+1}z)]\right]
 \label{bell}
\end{equation} 
for $(p, z) \in \mathbf{P} \times \mathbb{R}^{d}$,
acting on all functions $v$ where the stochastic kernel is
defined. Consider the \emph{Bellman recursion}, also referred to as
backward induction:
\begin{equation} 
v_{T}(p, z)= r_{T}(p,z), \quad v_{t}= {\cal T}_{t}
v_{t+1} \qquad \hbox{for $t=T-1, \dots,
0$.} \label{backward_induction}
\end{equation}
Having assumed that the reward functions are convex and
globally Lipschitz and the disturbances $W_{t+1}$ are integrable,
there exists a recursive solution $(v^{*}_{t})_{t=0}^{T}$ to the
Bellman recursion (\cite{hinz2014}). These functions
$(v^{*}_{t})_{t=0}^{T}$ are called \emph{value functions} and they
determine an optimal policy $\pi^{*}=(\pi^{*}_{t})_{t=0}^{T}$ via
\begin{align} 
\pi^{*}_{t}(p,z) & = \arg\max_{a \in
\mathbf{A}}\left[r_{t}(p, z, a)+ \sum_{p'\in\mathbf{P}}
\alpha_{p,p'}(a)\mathbb{E}[v^{*}_{t+1}(p',W_{t+1}z)] \right],
\end{align} 
for $t = T-1,\dots, 0$.

\section{Primal Approximate Solution  \label{sec-algo}}
The first step in obtaining a numerical solution to the backward
induction (\ref{backward_induction}) is an appropriate discretization
of the Bellman operator (\ref{bell}) to
\[{\cal T}^{n}_{t}v(p,z) = \max_{a \in \mathbf{A}}\left({r_{t}(p,z, a)+}
{\sum_{p'\in\mathbf{P}} \alpha_{p,p'}(a) \sum_{k=1}^{n}
\nu^n_{t+1}(k)v(p', W_{t+1}(k)z)}\right)\]
where the probability weights $(\nu^n_{t+1}(k))_{k=1}^{n}$ corresponds
to the distribution sampling $(W_{t+1}(k))_{k=1}^{n}$ of each
disturbance $W_{t+1}$. In the resulting modified backward induction
governed by $v^{n}_{t}={\cal T}^{n}_{t}v^{n}_{t+1}$, the modified
functions $(v^{n}_{t})_{t=0}^{T}$ need to be described by
algorithmically tractable objects. Since the reward and scrap
functions are convex in the continuous variable, these modified value
functions are also convex and can be approximated by piecewise linear
and convex functions. For this, introduce the so-called
sub-gradient envelope ${\cal S}_{\mathbf{G}^m}f$ of a convex function $f: \rr^d
\to \rr$ on a grid $\mathbf{G}^m \subset \rr^d$ with $m$ points 
i.e. $\mathbf{G}^{m}=\{g^{1}, \dots, g^{m}\}$
\[ {\cal S}_{\mathbf{G}^m}f=\vee_{g \in \mathbf{G}^m} (\triangledown_{g}f)\]
 which is a maximum of the sub-gradients $\triangledown_{g}f$ of $f$
on all grid points $g \in \mathbf{G}^m$. Using the sub-gradient envelope
operator, define the double-modified Bellman operator as
\[ {\cal T}^{m, n}_{t}v(p, z) = {{\cal S}_{\mathbf{G}^{m}}} {\max_{a \in
\mathbf{A}}}\left( {r_{t}(p, z , a)+}
{\sum_{p'\in\mathbf{P}} \alpha_{p,p'}(a)\sum_{k=1}^{n}\nu^n_{t+1}(k)v(p', W_{t+1}(k)z)} \right)
.\]
The corresponding backward induction
\begin{eqnarray} 
v^{m, n}_{T}(p, z)&=&{\cal S}_{\mathbf{G}^{m}} \max_{a\in
\mathbf{A}}r_{T}(p, z, a), \quad p \in \mathbf{P}, z \in \mathbb{R}^d \label{scheme1}\\
v^{m, n}_{t}(p, z)&=&{\cal T}^{m, n}_{t}v^{m, n}_{t+1}(p, z),
\qquad t=T-1, \dots 0. \label{scheme2}
\end{eqnarray}
yields the so-called double-modified value functions $(v^{n,
m}_{t})_{t=0}^{T}$. Under appropriate assumptions (see \cite{hinz2014}),
the double-modified value functions converge uniformly to the true
value functions almost surely on compact sets. These assumptions
include the convexity and global Lipschitz continuity of the rewards,
the integrability of all disturbances and some restrictions on the
distribution sampling and grid density.

Since the double-modified value functions $(v^{m, n}_{t})_{t=0}^{T}$
are piece-wise linear and convex, they can be expressed in a compact
and appealing form using matrix representations. Note that any
piecewise convex function $f$ can be described by a matrix where each
linear functional is represented by a row in the matrix. To denote
this relation, let us agree on the following notation: Given a
function $f$ and a matrix $F$, we write $f \sim F$ whenever $f(z)=\max
(Fz)$ holds for all $z\in\rr^d$. Let us emphasize that the
sub-gradient envelope operation ${\cal S}_{\mathbf{G}^m}$ on a grid
$\mathbf{G}^m$ is reflected in terms of a matrix representative by a
specific row-rearrangement operator
\[ f \sim F \quad \Leftrightarrow \quad {\cal S}_{\mathbf{G}^m}f \sim
\Upsilon_{\mathbf{G}^m}[F] \] 
where the row-rearrangement operator $\Upsilon_{\mathbf{G}^m}$
associated with $\mathbf{G}^m=\{g^{1}, \dots, g^{m}\} \subset \rr^{d}$ acts
on matrix $F$ with $d$ columns as follows:
\begin{equation} \label{row-rearrangment} 
  ({\Upsilon_{\mathbf{G}^m} }F)_{i, \cdot}= F_{ {\rm argmax}(F g^{i}),
    \cdot} \qquad \hbox{for all $i=1, \dots, m.$}
\end{equation}
For piecewise convex functions, the result of maximization, summation,
and composition with linear mapping, followed by sub-gradient envelope
can be obtained using their matrix representatives. More precisely, if
$f_{1} \sim F_{1}$ and $f_{2} \sim F_{2}$ holds, then it follows that
\begin{eqnarray*}
{\cal S}_{\mathbf{G}^m}(f_{1}+f_{2}) & \sim & \Upsilon_{\mathbf{G}^m}(F_{1})+ \Upsilon_{\mathbf{G}^m}(F_{2}) \\
{\cal S}_{\mathbf{G}^m}(f_{1} \vee f_{2}) & \sim & \Upsilon_{\mathbf{G}^m}(F_{1} \sqcup F_{2}) \\
{\cal S}_{\mathbf{G}^m}(f_{1}( W \cdot ) & \sim & \Upsilon_{\mathbf{G}^m}(F_{1}W)
\end{eqnarray*}
where $W$ is an arbitrary $d \times d$ matrix and the operator ${\sqcup}$ stands
for binding matrices by rows. Therefore, the backward induction
(\ref{scheme1}) and (\ref{scheme2}) can be expressed in terms of the
matrix representatives $V^{m, n}_{t}(p)$ of the value functions
$v^{(m, n)}_{t}(p, z)$ for $p \in \mathbf{P}$, $z \in \mathbb{R}^d$ and $t=0, \dots
T$. Since the double-modified backward induction involves
maximization, summations and compositions with linear mappings applied
to piecewise linear convex functions, it can be rewritten in terms of
matrix operations which results in the following algorithm.

\vspace{5mm}
\begin{minipage}{0.90\columnwidth}
\setlength{\algomargin}{-0.1em}
\begin{algorithm}[H]  
  \caption{Convex Stochastic Switching Algorithm\label{algo1}}
  \For {$p \in \mathbf{P}$}{
    $V_T^{m,n}(p)\leftarrow \Upsilon_{\mathbf{G}^m} \sqcup_{a\in \mathbf{A}} {\cal S}_{{\mathbf G}^m} r_{T}(p,.,a)$
  }
  \For{$t \in \{T-1,\dots,0\}$}{
    \For {$p \in \mathbf{P}$}{
      $\tilde V^{m,n}_{t+1}(p) \leftarrow \sum_{k=1}^n \nu^n_{t+1}(k) \Upsilon_{G^m} 
      V_{t+1}^{m,n}(p)  W_{t+1}(k)  $
    }
    \For {$p \in \mathbf{P}$}{
      $ V_{t}^{m,n}(p) \leftarrow \Upsilon_{\mathbf{G}^m} \sqcup_{a\in \mathbf{A}} \left( 
        {\cal S}_{{\mathbf G}^m} r_{t}(p,.,a) +
        \sum_{p'\in\mathbf{P}}\alpha_{p,p'}(a)\tilde V_{t+1}^{m,n}(p') \right)$
    }
  }
\end{algorithm}
\end{minipage}
\vspace{5mm}

The continuation value $\tilde V^{m,n}_{t+1}(p)$ given by  the sixth  line in
Algorithm \ref{algo1} can be further approximated by recalling that
the disturbances $W_{t+1}$ are independently and identically
distributed across time. Through the use of a suitable nearest
neighbour algorithm to construct permutation matrices, the conditional
expectation on the sixth line can be approximated to save a significant amount 
of computational effort. For further details, we refer the interested reader to 
\cite{hinz_yap2015}.

A candidate for a nearly optimal policy is given by
\[\pi^{m,n}_{t}(p,z)= \arg\max_{a \in \mathbf{A}}\left[r_{t}(p,z,a)
+ \max\left(\sum_{p'\in\mathbf{P}}\alpha_{p,p'}(a)\tilde V^{m,n}_{t+1}(p')z \right) \right].\] 
Given the pointwise convergence of the value functions for all $ p \in
\mathbf{P}$, $a\in \mathbf{A}$, and $t=0, \dots, T-1$ one can deduce that
\[\lim_{(m, n) \to \infty}\pi^{m,n}_{t}(p,z)=\pi^{*}(p,z) \quad p \in
\mathbf{P}, a\in \mathbf{A}, \enspace t=0, \dots, T-1\] 
holds for each point $(p, z) \in \mathbf{P} \times \rr^d$. However,
this convergence may be of little help in determining the quality of a
candidate policy in practice. To address this issue, a pathwise
dynamic approach will be described in the next subsection.  This area
has recently attracted growing research interest on diagnostics and
potential a posteriori justification of an approximate control policy.

\section{ Dual Solution And Diagnostics
  \label{sec-algo-diagnostic}}

This section utilizes an approach suggested in \cite{rogers2007} and
\cite{brown_etal2010} to assess the quality of a given control policy
obtained by primal methods.  Let us first illustrate their use for the
case of a optimal stopping problem. Given a real-valued stochastic
process $(Z_{t})_{t=0}^{T}$ adapted to a filtration, consider the
family ${\cal V}$ of all $\{0, 1, \dots, T\}$-valued stopping times.
Obviously, the optimal stopping value is dominated
\[V^{*}_{0}=\sup_{\tau \in {\cal V}} \EE(Z_{\tau}) \le
\EE(\sup_{0 \le t \le T} Z_{t})\]
by the expectation of the pathwise maximum. This dominance still holds
\begin{equation} 
V^{*}_{0} = \sup_{\tau \in {\cal V}}
\EE(Z_{\tau}-M_{\tau}) \le \EE( \sup_{0 \le t \le T} (Z_{t}-M_{t})).
\nonumber
\end{equation} 
when the original process $(Z_{t})_{t=0}^{T}$ is replaced by
$(Z_{t}-M_{t})_{t=0}^{T}$ with a martingale $(M_{t})_{t=0}^{T}$
starting at the origin $M_{0}=0$.  Furthermore, it turns out that this
estimate is tight
\begin{equation} 
  V^{*}_{0}=\inf_{(M_t)_{t=0}^T } \EE( \sup_{0
    \le t \le T}(Z_{t}-M_{t}))= \EE( \sup_{0
    \le t \le T}(Z_{t}-M^*_{t})) \label{martdom}
\end{equation}
where the infimum is taken over the family of all martingales starting
at the origin and in fact is attained at some {\it optimal martingale}
$(M^*_t)_{t=0}^T$ from this family. Usually, an optimal martingale is
not available, as its knowledge is equivalent to the solution of the
stopping problem.  The equation (\ref{martdom}) yields a practical way
to estimate the optimal value $V^{*}_{0}$ via a Monte Carlo procedure
by determining an average of the pathwise maximum on a number of
simulated trajectories of $(Z_{t}- \tilde M_{t})_{t=0}^{T}$.
Obviously, to obtain a tight upper estimate, the martingale $(\tilde
M_{t})_{t=0}^{T}$ must resemble the optimal one.  Similarly,
a lower estimate for $ V^{*}_{0}$ can be obtained from Monte-Carlo
average of independent realizations of $Z_{\tilde \tau}- \tilde
M_{\tilde \tau}$, based on an arbitrary stopping time $\tilde \tau$.
Again to obtain a tight lower bound, the stopping time $\tilde\tau$
must be chosen close to the optimal stopping time.  This procedure of
bound estimate exhibits an interesting self-tuning property: The
closer $(\tilde M_{t})_{t=0}^{T}$ and $\tilde \tau$ are to their
optimal counterparts, the lower the variance of the Monte-Carlo
trials and the narrower the gap between both bounds.  In a
hypothetic case where $\tilde \tau$ is the optimal stopping time and
$(\tilde M_{t})_{t=0}^{T}$ is the optimal martingale, the variance of
both Monte-Carlo tails reduces to zero. Moreover, in this case the
upper and the lower bounds coincide with the optimal stopping value.
 
This approach has been successfully applied
\cite{andersen_broadie2004} and generalized to multiple stopping
problems in finance \cite{meinshausen_hambly2004}. For our switching
problem, the pathwise dynamic programming approach works similarly
(see \cite{hinz_yap2015}).  Given an approximate control problem solution,
two random bound variables are constructed, whose expectations give a
lower and an upper estimate for the unknown value function.  These
bound variables are determined through a stochastic recursive
procedure similar to the backward induction. Their calculation is
effected in terms of a Monte-Carlo procedure whose in-sample empirical
variance can be used to quote a confidence interval for the unknown
value function.  In this setting, the self-tuning property states that
both, the variance of the Monte-Carlo trials and the gap between
expectations of bound variables decrease, if the {\it approximate
  control policy } approaches the optimal policy.

Suppose that our numerical scheme returns approximate value
functions $(v_{t})_{t=0}^{T}$ along with the approximate expected
value functions $(v^{E}_{t})_{t=0}^{T}$, thus we introduce an
approximately-optimal policy $\pi=(\pi_{t})_{t=0}^{T-1}$ by
	$$
	\pi_{t}(p,z)={\rm argmax}(r_{t}(p,z, a)+ \sum_{p'\in
\mathbf{P}}\alpha_{p, p'}(a)v^{E}_{t+1}(p', z)))
	$$
   To answer the question how far is the strategy $\pi$ from an
optimal strategy $\pi^*$, we estimate the performance gap $ {[v^{
\pi}_{0}(p_{0}, z_{0}), v^{\pi^{*}}_{0}(p_{0}, z_{0})]} $ at a given
point a point $(p_{0}, z_{0}) \in \mathbf{P}\times \rr^d$ by an
explicit construction of random variables $ \underline{v}^{\pi,
\varphi}_{0}(p_{0}, z_{0})$, $\overline{v}^{ \varphi}_{0}(p_{0},
z_{0}) $ satisfying
	$$
	\EE(\textcolor[rgb]{0.00,0.0,0.00}{\underline v^{\pi,
\varphi}_{0}(p_{0}, z_{0})})=\textcolor[rgb]{0.00,0.00,0.00}{v^{
\pi}_{0}(p_{0}, z_{0}) \le v^{\pi^{*}}_{0}(p_{0}, z_{0})} \le
\EE(\textcolor[rgb]{0.00,0.0,0.00}{\bar v^{\varphi}_{0}(p_{0},
z_{0})}).
	$$
	Using Monte-Carlo simulations, one estimates both means along
with empirical confidence intervals to estimate the performance gap,
which is narrow if the approximate solution is close to optimal. This
useful property is due to the the self-tuning property, which ensures
that the more optimal solutions $(v_{t})_{t=0}^{T}$
$(v^{E}_{t})_{t=0}^{T}$ return narrower gaps and lower variance in the
Monte-Carlo scheme. This scheme can be implemented as follows
	
\begin{enumerate}
\item Given approximate solution $(v_{t})_{t=0}^{T}$  $(v^{E}_{t})_{t=0}^{T}$
with the corresponding policy $( \pi_{t})_{t=0}^{T-1}$, implement control
variables  $(\varphi_{t})_{t=1}^{T}$ as
$$
\varphi_{t}(p, z, a)=\sum_{p'\in \mathbf{P}}\alpha_{p,p'}(a)\left(\frac{1}{I}
\sum_{i=1}^{I} v_{t}(p', W^{(i)}_{t}z) - v_{t}(p', W_{t}z)\right),
$$
for all $p \in \mathbf{P}$, $a \in \mathbf{A}$, $z \in \rr^{d}$, where $(W^{(1)}_{t},
\dots, W^{(I)}_{t}, W_{t})_{t=1}^T$ are independent and $(W^{(1)}_{t},
\dots, W^{(I)}_{t}, W_{t})$ are identically distributed for each $t=1,
\dots, T$.
\item Chose a number $K \in \nn$ of Monte-Carlo trials and obtain for
$k=1, \dots, K$ independent realizations $
(W_{t}(\omega_{k}))_{t=1}^{T} $ of disturbances.
\item Starting at $z_{0}^{k}:=z_{0} \in \rr^{d}$, define for $k=1,
\dots, K$ trajectories $(z_{t}^{k})_{t=0}^{T}$ recursively
$$
z^{k}_{t+1}=W_{t+1}(\omega_{k})z^{k}_{t}, \qquad t=0, \dots, T-1
$$
and determine realizations
$$
\varphi_{t}(p, z^{k}_{t-1}, a)(\omega_{k}), \qquad t=1, \dots, T,\quad k=1, \dots, K.
$$
\item For each $k=1, \dots, K$ initialize the  recursion at $t=T$ as
  $$
  \underline{v}^{ \pi,  \varphi}_{T}(p,z^{k}_{T})(\omega_{k})=
  \overline{v}^{ \pi}_{T}(p,z^{k}_{T})(\omega_{k})=
  r_{T}(p, z^{k}_{T}) \qquad \hbox{for all } p \in \mathbf{P}
  $$
  and continue for $t=T-1, \dots, 0$ and for all $p \in \mathbf{P}$ by
  \begin{eqnarray*} 
    \lefteqn{
    \underline{v}^{ \pi,  \varphi}_{t}(p,z^{k}_{t})(\omega_{k})= r_{t}(p, z^{k}_{t}, \pi_{t}(p, z^{k}_{t}))+\varphi_{t+1}(p,z^{k}_{t},  \pi_{t}(p, z^{k}_{t}))(\omega_{k})}  \qquad \qquad \qquad\\
 && \qquad \qquad \qquad
    + \sum_{p' \in \mathbf{P}} \alpha_{p,p'}(\pi_{t}(p, z^{k}_{t})) \underline{v}^{ \pi, \varphi}_{t+1}(p', z^{k}_{t+1})(\omega_{k})
    \nonumber \label{rec2-1} \\
    \lefteqn{
    \overline{v}^{\varphi}_{t}(p,z^{k}_{t})(\omega_{k})= \max_{a \in \mathbf{A}} \big [ r_{t}(p, z^{k}_{t}, a) + \varphi_{t+1}(p,z^{k}_{t},a)(\omega_k)}  \qquad \qquad \qquad\\
 && \qquad \qquad \qquad
    + \sum_{p' \in \mathbf{P}} \alpha_{p,p'}(a) \overline{v}^{\varphi}_{t+1}(p', z^{k}_{t+1})(\omega_{k}) \big ]
    \nonumber \label{rec2-1}
  \end{eqnarray*}
\item Calculate sample means
  $\frac{1}{K}\sum_{k=1}^{K} \underline{v}^{\pi, \varphi}_{0}(p_{0},
  z_{0})(\omega_{k})$, $\frac{1}{K}\sum_{k=1}^{K}
  \overline{v}^{\varphi}_{0}(p_{0}, z_{0})(\omega_{k})$ to estimate the
  means $\EE(\underline{v}^{\pi, \varphi}_{0}(p_{0}, z_{0}))$,
  $\EE(\overline{v}^{\varphi}_{0}(p_{0}, z_{0}))$ along with their
  confidence intervals. These means will be referred to as the \emph{primal}
  and \emph{dual} values, respectively.
\end{enumerate}

\section{Numerical Results \label{sec-result}}
This section demonstrates our algorithms on the commodity extraction
problem described in Section \ref{sec-problem}. The values given by
Table \ref{bs_table} will be adjusted in order fit our discrete time
setting with time step $\Delta = 0.25$ years between the decision
times $\{ 0,\dots ,T\} \subset \nn$ representing a time horizon $[0,
\bar T] = [0, 30]$ of thirty years, with four decisions made every
year.  Let us discretize the commodity levels remaining in the
resource with the step size corresponding to commodity amount
extracted within one time step $\Delta$.  That is, the first entry of
the discrete component covers the range $p^{(1)} \in \left\{0, 1, ...,
\frac{\bar p}{\Delta}\right\}$ where ${\bar p}=15$ stands for the
minimum number of years it takes to deplete the resource. Further,
define the number of decision epochs $T = \frac{\bar T}{\Delta} + 1$,
the maintenance costs $m_{t} = m_0 \Delta e^{(\rho - r - \zeta)
t\Delta}$ and the switching costs $c_{t} = c_0 e^{(\rho - r - \zeta)
t\Delta}$ with coefficients $r = 0.1$, $\rho = 0.08$, $\zeta = 0.02$,
$m_0 = 0.5$ and $c_0 = 0.2$ standing for the interest rate, rate of
inflation, real estate tax, initial maintenance cost and initial 
switching cost, respectively.  \\

\noindent{\bf Geometric Brownian Motion:}
Following \cite{brennan_schwartz1985}, the continuous
state component $(Z_t)_{t=0}^T$ is one-dimensional and follows
linear state dynamics:
\begin{equation} \label{bsdynamics}
  Z_{t+1} = \exp\left(\left(\mu - \frac{\sigma^2}{2}\right)\Delta +
      \sigma\sqrt{\Delta} N_{t+1} \right) Z_{t}
  \end{equation}
  with $\sigma^2 = 0.08$ and independent identically standard normally
distributed $(N_t)_{t=1}^T$. The controlled transition probabilities of
the discrete component given $(p^{(1)}, p^{(2)}) \in
\mathbf{P}$ is uniquely determined by
  \begin{align*}
    & \alpha_{(p^{(1)},p^{(2)}),(\max\{p^{(1)}-1,0\},\text{Opened})}(\text{Open}) = 1,\\
    & \alpha_{(p^{(1)},p^{(2)}),(p^{(1)},\text{Closed})}(\text{Close}) = 1,\\
    & \alpha_{(p^{(1)},p^{(2)}),(0,p^{(2)})}(\text{Abandon}) = 1.
  \end{align*}
The mode is controlled deterministically under the above specfication.
Now let us define the control costs by a function, representing
the revenue of the opened resource
\[h_{t}(z) = 5\Delta ze^{-(r+\zeta)t\Delta} -2.5\Delta e^{(\rho - r - \zeta)t\Delta}\]
where the variable $z$ represents the commodity price.

An equally spaced grid ranging from 0 to 20 of size 4,001 was used and
a distribution sampling of size 20,000 constructed using equidistant
sampling of standard normal quantiles was applied to obtain the
following results.  Table \ref{table_value1} lists the primal and dual
values, generated using $K=1,000$ paths and $I=1,000$ sub-simulations.
To compare our results with \cite{brennan_schwartz1985} listed in the
first column of the Table \ref{table_value1}, we assumed a convenience
yield $d = 0.01$ and set $\mu = r - d = 0.09$ in
(\ref{bsdynamics}). We were able to obtain similar results despite
\cite{brennan_schwartz1985} using a continuous time formulation.
Further, our results were able to obtain excellent precision
as evidenced by the tight bounds and low standard errors given in the
paranthesis.

\begin{table}[h]
  \caption{Resource Valuation Under Geometric Brownian Motion 
    \label{table_value1}}
  \centering
  \setlength\tabcolsep{2pt} 
  \begin{tabular}{l|c|c|c|c|c|c}
    &\multicolumn{2}{c}{B\&S}&\multicolumn{2}{c}{Opened Resource}&\multicolumn{2}{c}{Closed Resource}\\
    $Z_0$ & Open & Closed & Primal & Dual & Primal & Dual \\
    \hline
    0.3 & 1.25  & 1.45  &  1.2127(.0026) &  1.2156(.0026)&  1.4127(.0026) &  1.4156(.0026) \\
    0.4 & 4.15  & 4.35  &  4.1059(.0034) &  4.1086(.0034)&  4.3059(.0034) &  4.3086(.0034) \\
    0.5 & 7.95  & 8.11  &  7.9026(.0039) &  7.9053(.0039)&  8.0752(.0041) &  8.0777(.0041) \\
    0.6 & 12.52 & 12.49 & 12.5129(.0042) & 12.5153(.0042)& 12.4787(.0045) & 12.4813(.0045) \\
    0.7 & 17.56 & 17.38 & 17.5869(.0047) & 17.5889(.0047)& 17.3869(.0047) & 17.3904(.0048) \\
    0.8 & 22.88 & 22.68 & 22.9475(.0052) & 22.9489(.0052)& 22.7475(.0052) & 22.7489(.0052) \\
    0.9 & 28.38 & 28.18 & 28.4940(.0057) & 28.4957(.0057)& 28.2940(.0057) & 28.2957(.0057) \\	
    1.0 & 34.01 & 33.81 & 34.1667(.0062) & 34.1681(.0062)& 33.9667(.0062) & 33.9681(.0062) \\
  \end{tabular}
\end{table}

\begin{figure}[h]
  \centering
  \includegraphics[height=2.3in,width=2.3in]{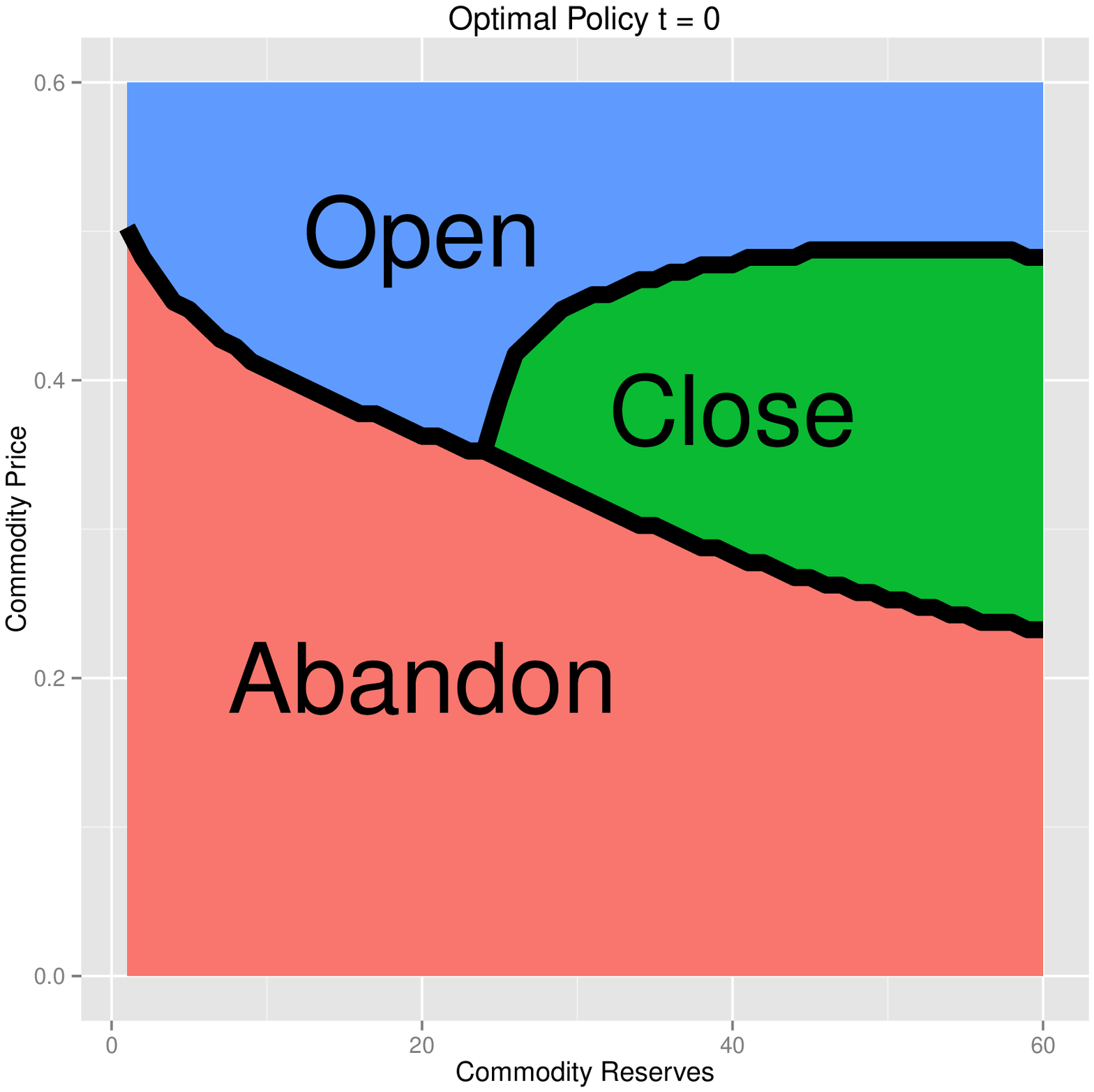}
  \includegraphics[height=2.3in,width=2.3in]{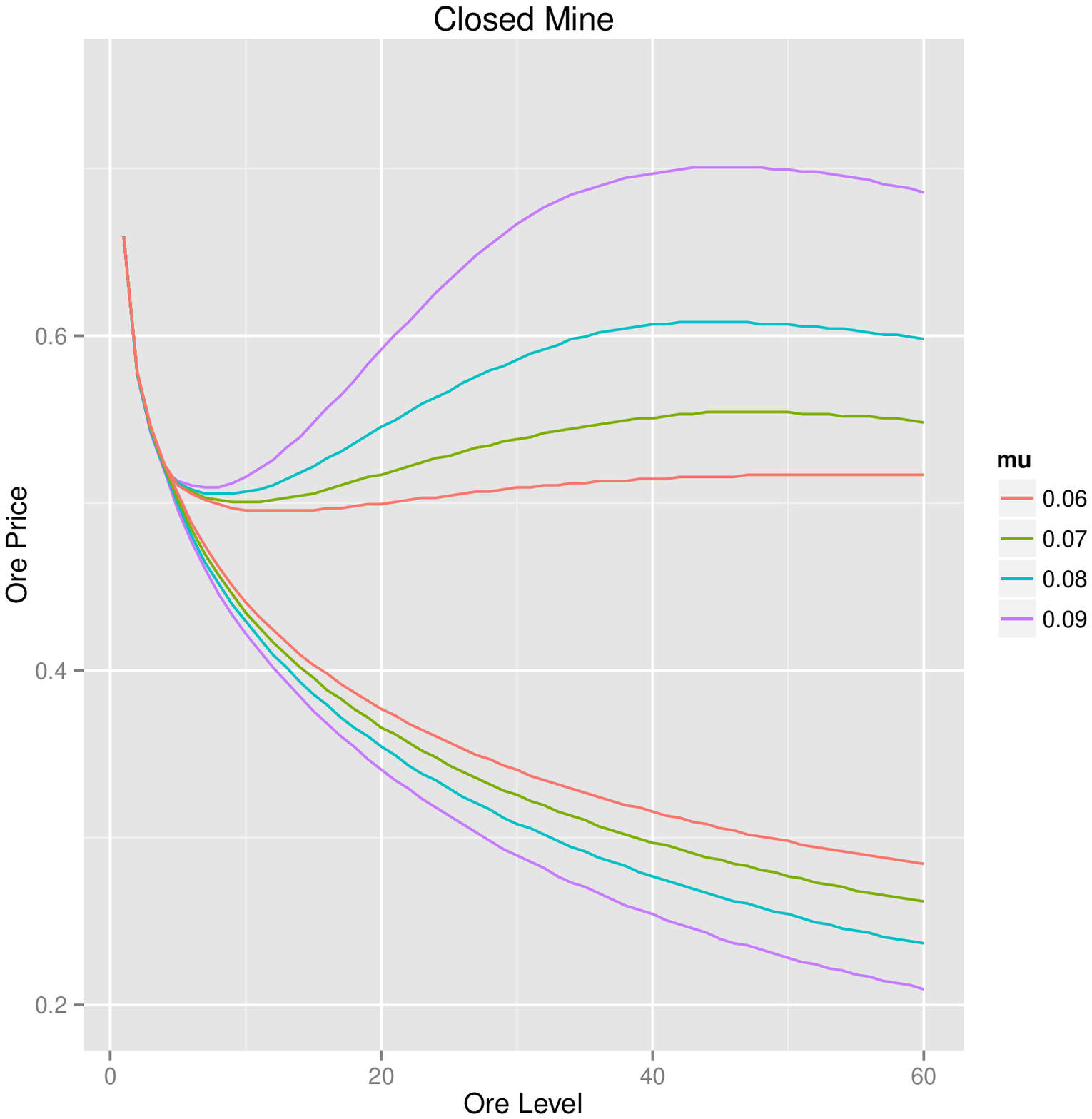}
  \caption{The optimal policy at $t = 0$. The horizontal axis represents
    remaining reserves while the vertical represents the commodity price.
    The right plot displays policy for a commodity resource that begins closed for
    $\mu = 0.09, 0.08, 0.07, 0.06$.}
  \label{plot1}
\end{figure}

This paper will illustrate the control policies in the form of the
curves depicted in Figure \ref{plot1}. We have observed that for all
the situations examined in this paper, the optimal policies share this
structue. Namely, it is optimal to operate an open resource when the
commodity price is sufficiently high.  Below some price level, the
resource should be shut down, either to the closed or abandoned
mode. Typically, we have a curve bifurcation and the area between
both curve branches represents the price range where switching to the
closed mode is optimal. In the right plot of Figure \ref{plot1}, the
optimal policies for a closed resource is shown for different values
of $\mu$. Recall that there is a strong element of scarcity in our
problem. That is, selling a unit of limited commodity now means that
the same unit of commodity cannot be sold at a future date. Thus, we observe from
Figure \ref{plot1}, that the decision maker is more willing to
preserve commodity for future use if the drift term $\mu$ is high.
In our subsequent studies, we will assume drift $\mu = r - d = 0.09$ 
(unless otherwise stated) in order to stay in line with 
\cite{brennan_schwartz1985}.
\\

\begin{figure}[h!] 
  \centering
  \includegraphics[height=2.3in,width=2.3in]{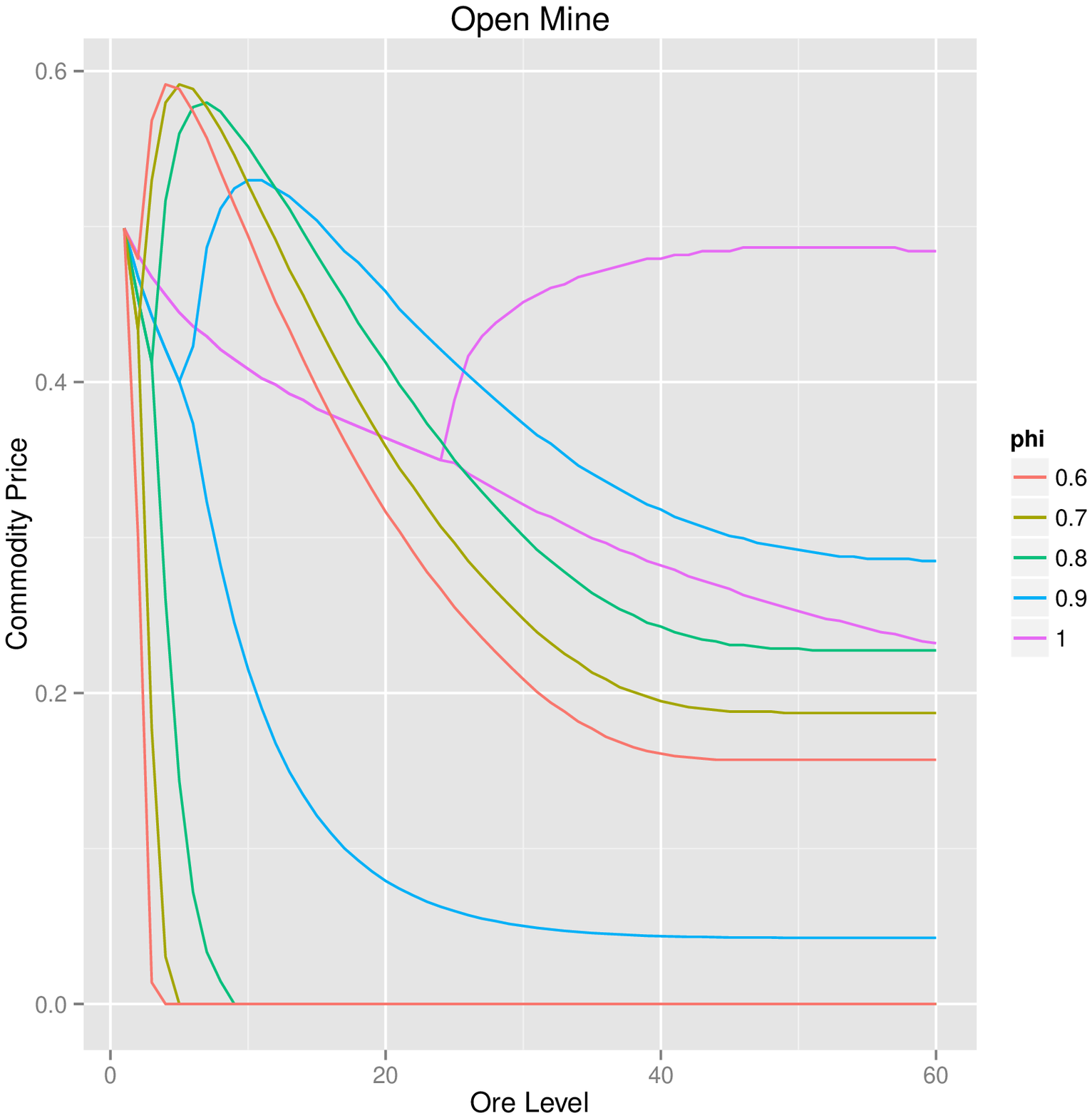}
  \includegraphics[height=2.3in,width=2.3in]{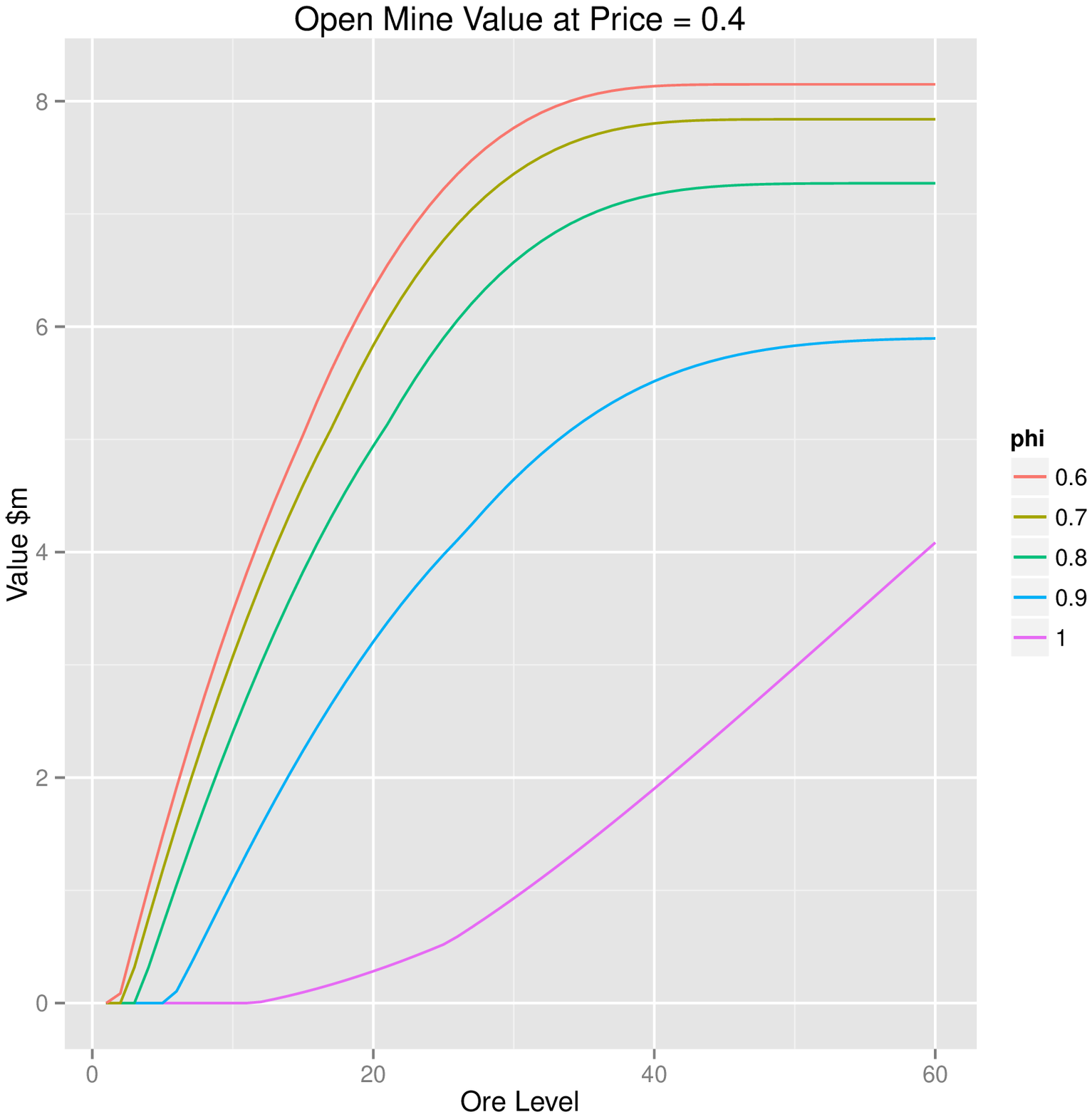}
  \caption{Left plot shows policy for an opened commodity resource under
    $\phi = 1, 0.9, 0.8, 0.7, 0.6$ while the right displays the value of
    the resource as a function of commodity reserves for $\tilde Z_0 = 0.4$.
}
  \label{plot_ar}
\end{figure}

\noindent
{\bf Ornstein-Uhlenbeck Process:}
Commodity prices often mimic the business cycle in practice (see
\cite{paschke_prokopczuk2010, schwartz1997}). Prices with significant
mean reversion often lead to different optimal behavior
compared to those following geometric Brownian motions 
(\cite{sarkar2003, schwartz1997, tsekrekos2010}). To examine the
impact of mean reversion, the logarithm of the commodity price $(\log
\tilde Z_t)_{t=0}^T$ is assumed to follow an AR(1) process, an
auto-regression of order one. With this choice, define the linear
state evolution $( Z_t)_{t=0}^T$ in the required form
$Z_{t+1}=W_{t+1}Z_t$ as
  \[Z_{t+1} := \begin{bmatrix} 1 \\ \log \tilde Z_{t+1} \end{bmatrix}
    = \begin{bmatrix} 1 & 0 \\ 
      \left(\mu - \frac{\sigma^2}{2} \right)\Delta + \sigma\sqrt{\Delta}N_{t+1} & \phi \end{bmatrix}
    \begin{bmatrix} 1 \\ \log \tilde Z_t \end{bmatrix},\]
for $ t=0, \dots, T-1,$ with initial value $ \log \tilde Z_{0} \in
\mathbb{R}$ and $\mu = 0.09$, $\sigma^2 = 0.08$. Again,
$(N_t)_{t=1}^T$ are independent standard normally distributed random
variables and the parameter $\phi \in [0, 1]$ determines the speed of
mean reversion, where $\phi = 1$ gives a geometric Brownian motion.
Using the same parameters as before, merely adjust the cash flow
function accordingly for $t=0, \dots, T-1$ to
\[h_{t}(z) = 5\Delta \exp\left(z^{(2)}\right)\exp({-(r+\zeta)t\Delta})
-2.5\Delta \exp((\rho - r - \zeta)t\Delta)\] since now the variable
$z^{(2)}$ captures the logarithmic commodity price.  Figure
\ref{plot_ar} shows the nearly-optimal policies at initial time $t=0$
for an opened asset, which are computed using 2000 equally spaced grid
points ranging from -5 to 5 for $\log (\tilde Z_t)$ and 10,000
disturbances constructed from the equidistant quantiles of the normal
distribution.  Diagnostics is performed using $K=500$ sample paths and
$I= 500$ sub-simulations. Table \ref{table_ar} contains results for
$\phi = 1, 0.8, 0.6$ and shows very tight bounds and low standard
errors.

\begin{table}[h!]
  \centering
  \setlength\tabcolsep{8pt}
  \caption{Mean Reversion On Commodity Resource Valuation \label{table_ar}}
  \begin{tabular}{lccccc}
    & &\multicolumn{2}{c}{Opened Resource}&\multicolumn{2}{c}{Closed Resource}\\
    $\phi$ & $\tilde Z_0$ & Primal & Dual & Primal & Dual \\
    \hline
      1 & 0.3 &  1.2198(.0047) &  1.2278(.0049)&  1.4198(.0049) &  1.4279(.0049) \\
        & 0.4 &  4.1153(.0067) &  4.1203(.0067)&  4.3153(.0067) &  4.3203(.0067) \\
        & 0.5 &  7.9133(.0074) &  7.9179(.0075)&  8.0863(.0078) &  8.0900(.0078) \\
    \hline
    0.8 & 0.3 &  6.3187(.0006) &  6.3189(.0006)&  6.3027(.0006) &  6.3029(.0006) \\
        & 0.4 &  7.2752(.0006) &  7.2753(.0006)&  7.0822(.0006) &  7.0823(.0006) \\
        & 0.5 &  8.1160(.0006) &  8.1161(.0006)&  7.9160(.0006) &  7.9161(.0006) \\
    \hline
    0.6 & 0.3 &  7.6732(.0003) &  7.6733(.0003)&  7.6034(.0003) &  7.6034(.0003) \\
        & 0.4 &  8.1514(.0003) &  8.1515(.0003)&  7.9535(.0003) &  7.9536(.0003) \\
        & 0.5 &  8.5754(.0003) &  8.5755(.0003)&  8.3754(.0003) &  8.3755(.0003) \\
  \end{tabular}
\end{table}

Figure \ref{plot_ar} provides an interesting insight. We observe that
under mean reversion $|\phi|<1$, each additional unit of commodity
reserve exhibit diminishing marginal value which directly
contrasts the behaviour under geometric Brownian motion. This
phenomenon was also observed for the closed commodity resource and
will have significant implications later on.
\\

\textbf{Wasteful Extraction:}
For $w \in [0,1]$ and commodity level $p^{(1)}$, assume that transition
to the new positon is now governed by the probabilities
\begin{align*}
& \alpha_{(p^{(1)},p^{(2)}),(\max\{p^{(1)}-2,0\},\text{Opened})}(\text{Open}) = w\\
& \alpha_{(p^{(1)},p^{(2)}),(\max\{p^{(1)}-1,0\},\text{Opened})}(\text{Open}) = 1 - w\\
& \alpha_{(p^{(1)},p^{(2)}),(p^{(1)},\text{Closed})}(\text{Close}) = 1\\
& \alpha_{(p^{(1)},p^{(2)}),(0,p^{(2)})}(\text{Abandon}) = 1
\end{align*}
and zero otherwise, i.e. there is a probability $w$ of wasting 1
additional unit of the commodity when the resource is in the opened
mode.  

\begin{table}[h!]
	\caption{Commodity Resource Valuation Under Wastage \label{table_waste}}
	\begin{tabular}{lcccccc}
		& & &\multicolumn{2}{c}{Opened Resource}&\multicolumn{2}{c}{Closed Resource}\\
		$\phi$ & $w$ & $\tilde Z_0$ & Primal & Dual & Primal & Dual \\
		\hline
		1  &    0 &   0.3 &  1.2070(.0052) & 1.2129(.0052) & 1.4070(.0052) & 1.4129(.0052) \\
		&      &   0.4 &  4.0989(.0071) & 4.1020(.0071) & 4.2989(.0071) & 4.3020(.0071) \\
		&      &   0.5 &  7.8935(.0081) & 7.9015(.0082) & 8.0660(.0085) & 8.0720(.0085) \\
		\cline{2-7}
		&  0.5 &   0.3 &  0.1987(.0032) & 0.2042(.0030) & 0.3987(.0032) & 0.4040(.0030) \\
		&      &   0.4 &  1.9322(.0047) & 1.9352(.0048) & 2.1322(.0047) & 2.1352(.0048) \\
		&      &   0.5 &  4.5048(.0051) & 4.5094(.0051) & 4.6702(.0056) & 4.6741(.0055) \\
		\hline
		0.6  &    0 &   0.3 &  7.6727(.0003) & 7.6728(.0003) & 7.6028(.0003) & 7.6029(.0003) \\
		&      &   0.4 &  8.1509(.0003) & 8.1510(.0003) & 7.9530(.0003) & 7.9531(.0003) \\
		&      &   0.5 &  8.5749(.0004) & 8.5750(.0004) & 8.3749(.0004) & 8.3750(.0004) \\
		\cline{2-7}
		&  0.5 &   0.3 &  7.6514(.0003) & 7.6515(.0003) & 7.5897(.0003) & 7.5897(.0003) \\
		&      &   0.4 &  8.1296(.0003) & 8.1296(.0003) & 7.9393(.0003) & 7.9393(.0003) \\
		&      &   0.5 &  8.5536(.0003) & 8.5536(.0003) & 8.3536(.0003) & 8.3536(.0003)
	\end{tabular}
\end{table}

\begin{figure}[h!]
	\includegraphics[height=2.3in,width=2.5in]{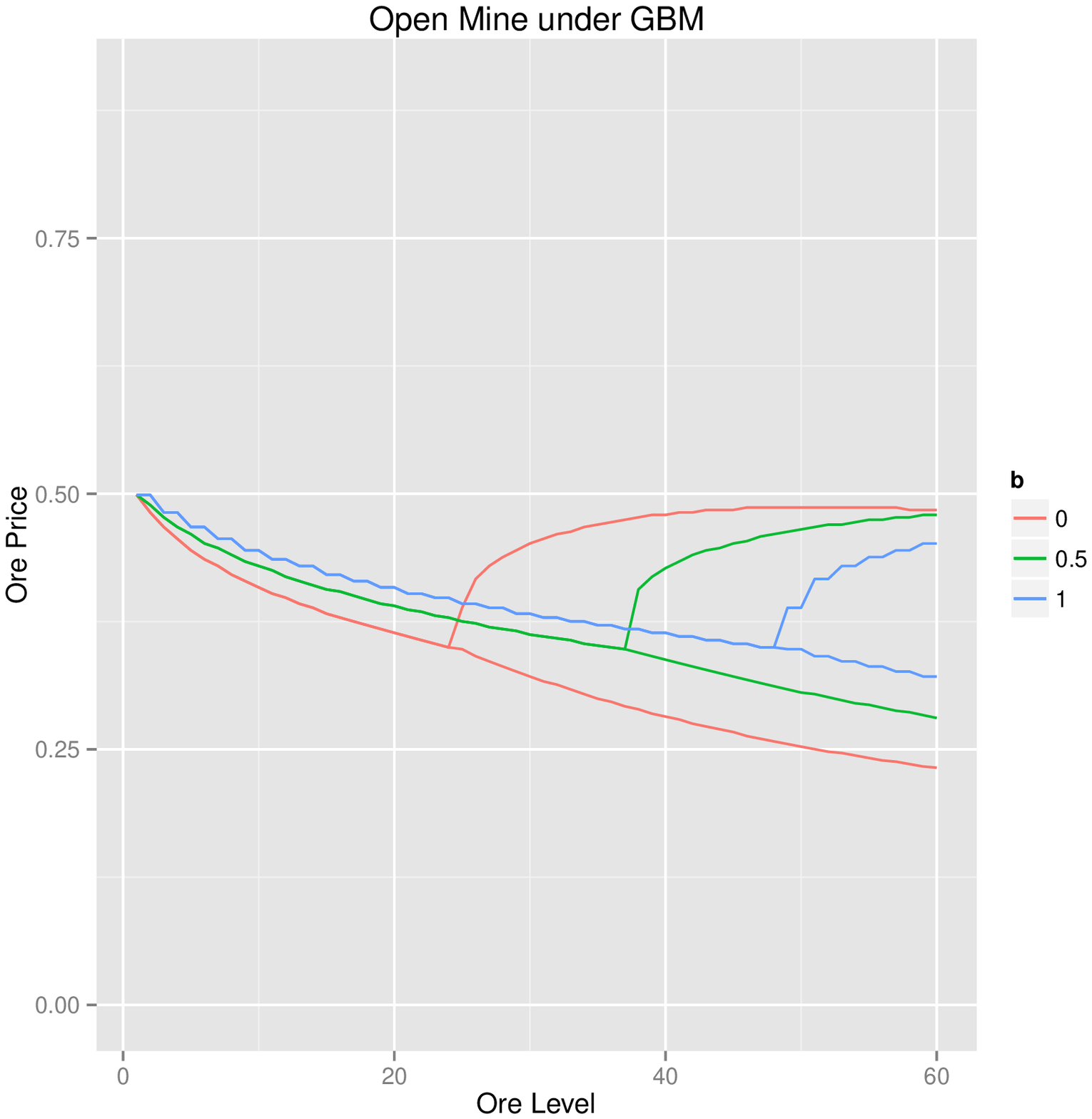}
	\includegraphics[height=2.3in,width=2.5in]{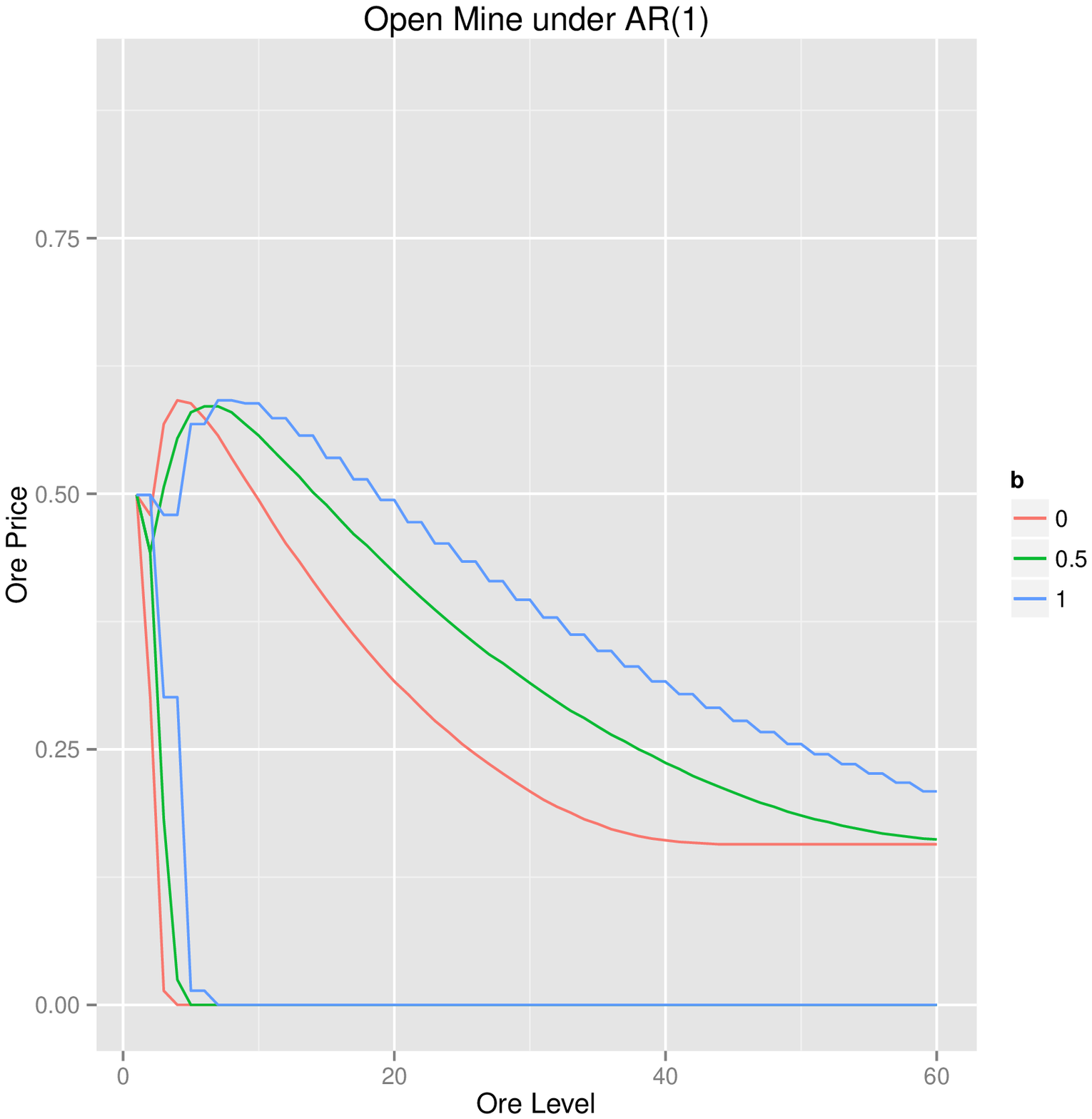}
	\caption{Optimal policies at $t = 0$ for an opened resource under
		$w = 0, 0.5, 1$. The policies under $\phi=1$ are depicted on the
		left while the right shows $\phi=0.6$. \label{plot_waste}}
\end{figure}

The linear state evolution for $(Z_t)_{t=0}^T$, the cash flow $h_t$,
the grid and the disturbance sampling of the disturbances from the
Ornstein-Uhlenbeck process study will be reused to generate the
results in this part. The primal and dual values in Table
\ref{table_waste} were generated using $K = 500$ sample paths and
$I = 500$ sub-simulations. As usual, the bounds obtained are tight
with low standard errors.  Table \ref{table_waste} suggests that
wasteful extraction has a minimal impact on the resource value when
the price follows strong mean reversion ($\phi = 0.6$) compared to
geometric Brownian motion for our prices starting at $\tilde Z_0$.
This may be directly related to the phenomenon observed earlier
regarding the marginal value of each additional unit of the
commodity. The impact on the optimal policies are shown in Figure
\ref{plot_waste} for $w = 0, 0.5, 1$.
\\

\textbf{Physical Delivery:} Another prominent feature of commodity
markets is the use of physical delivery contracts between suppliers
and buyers.  The reward function is modified in the following manner
to incorporate this. If the asset is abandoned, there are
neither costs nor revenue
\[ r_{t}((0, p^{(2)}, z), a)=0, \qquad a \in \mathbf{A}, z \in \mathbb{R} \]
For the case where $p^{(1)} > 0$, we have
\[r_{t}((p^{(1)}, p^{(2)}, z), a) = h_{t}(z)1_{\{2\}} (a) + m_{t} 1_{
	\{1\}} (a) + c_{t}1_{\{1,2\}}(a)|p^{(2)}-a| - \psi_t(p^{(1)},z)\]
where the values of $m_t, c_t$ are already specified before and
$\psi_t$ represents the penalty function for deviating from the
contract. This paper will consider penalty functions $\psi: \mathbf{P}
\times \mathbb{R}^2 \rightarrow \mathbb{R}$ of the form
\[\psi_t(p,z) = 
\begin{cases} 
b\exp(z^{(2)})(p^{(1)} - p_t^*), & \mbox{if } p^{(1)} > p_t^*;\\
0 & \mbox{otherwise} \end{cases}\]
where $z^{(2)}$ is the log commodity price, $b \geq 0$ and $p_0^{(1)}-p_t^*$
is the minimum amount of commodity that needs to be delivered by time $t$
according to the contract. If this condition is not feasible, then  the
decision maker will need to compensate the buyer with a money amount
equal to a proportion $b$ of the current market value of the
shortfall.

\begin{table}[h!]
	\setlength\tabcolsep{4pt} 
	\caption{Resource Valuation Under Physical Delivery Contracts 
		\label{table_contract}}
	\begin{tabular}{lcccccc}
		& & &\multicolumn{2}{c}{Opened Resource}&\multicolumn{2}{c}{Closed Resource}\\
		$\phi$ & $b$ & $\tilde Z_0$ & Primal & Dual & Primal & Dual \\
		\hline
		1   &    0 &   0.3 &  1.2070(.0052) & 1.2129(.0052) &  1.4070(.0052) & 1.4129(.0052) \\
		&      &   0.4 &  4.0989(.0071) & 4.1020(.0071) &  4.2989(.0071) & 4.3020(.0071) \\
		&      &   0.5 &  7.8935(.0081) & 7.9015(.0082) &  8.0660(.0085) & 8.0720(.0085) \\
		\cline{2-7}
		&    1 &   0.3 &  0.3133(.0036) & 0.3169(.0035) &  0.3703(.0036) & 0.3752(.0035) \\
		&      &   0.4 &  3.0951(.0056) & 3.1004(.0056) &  2.9738(.0057) & 2.9791(.0057) \\
		&      &   0.5 &  7.2354(.0070) & 7.2395(.0070) &  7.0435(.0071) & 7.0496(.0071) \\
		\hline
		0.6  &    0 &   0.3 &  7.6727(.0003) & 7.6728(.0003) &  7.6028(.0003) & 7.6029(.0003) \\
		&      &   0.4 &  8.1509(.0003) & 8.1510(.0003) &  7.9530(.0003) & 7.9531(.0003) \\
		&      &   0.5 &  8.5749(.0004) & 8.5750(.0004) &  8.3749(.0004) & 8.3750(.0004) \\
		\cline{2-7}
		&    1 &   0.3 &  7.6727(.0003) & 7.6728(.0003) &  7.5997(.0003) & 7.5998(.0003) \\
		&      &   0.4 &  8.1509(.0003) & 8.1510(.0003) &  7.9529(.0003) & 7.9530(.0003) \\
		&      &   0.5 &  8.5749(.0004) & 8.5750(.0004) &  8.3749(.0004) & 8.3750(.0004)
	\end{tabular}
\end{table} 

\begin{figure}[h!]
	\includegraphics[height=2.3in,width=0.5\textwidth]{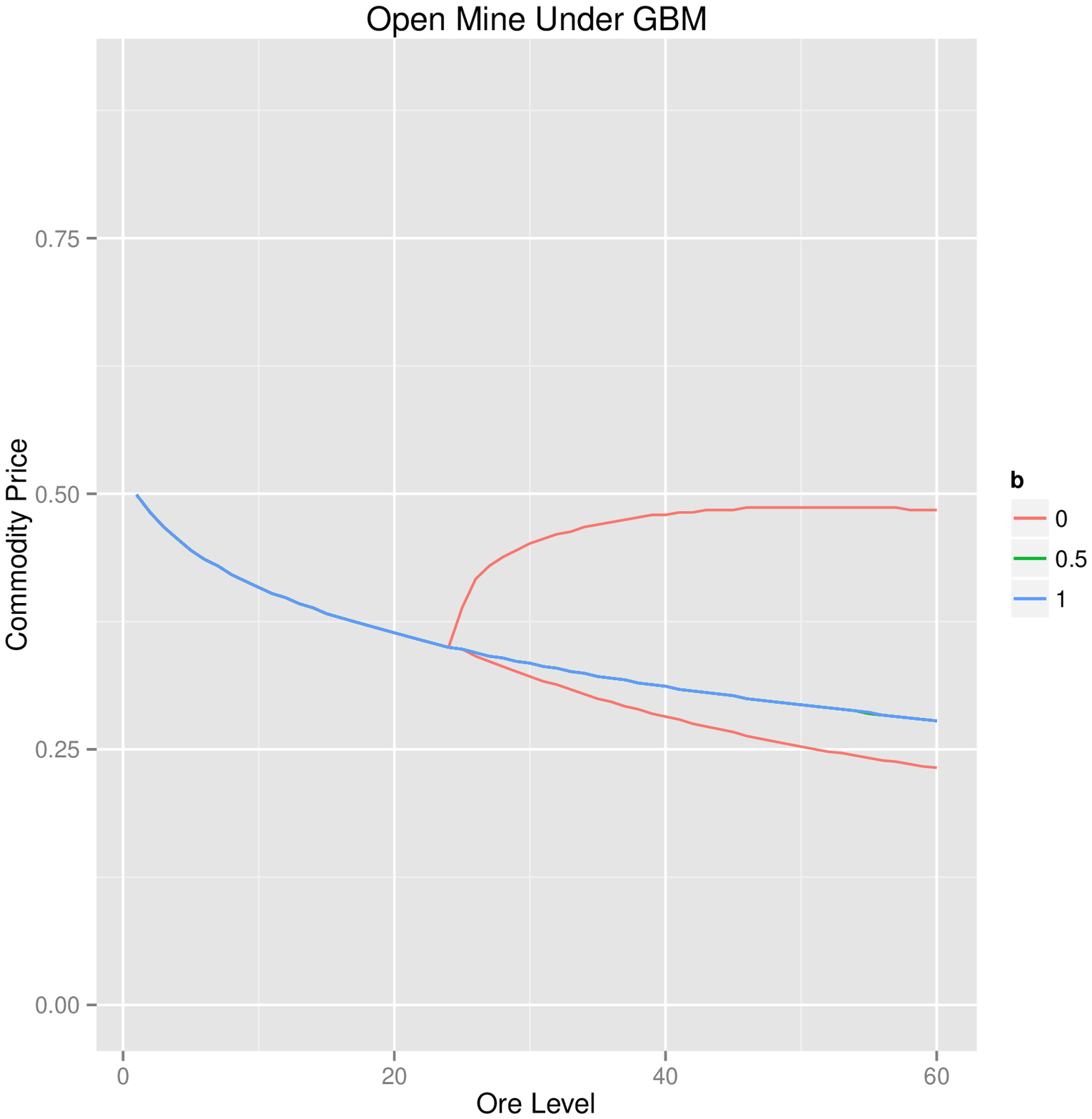}
	\includegraphics[height=2.3in,width=0.5\textwidth]{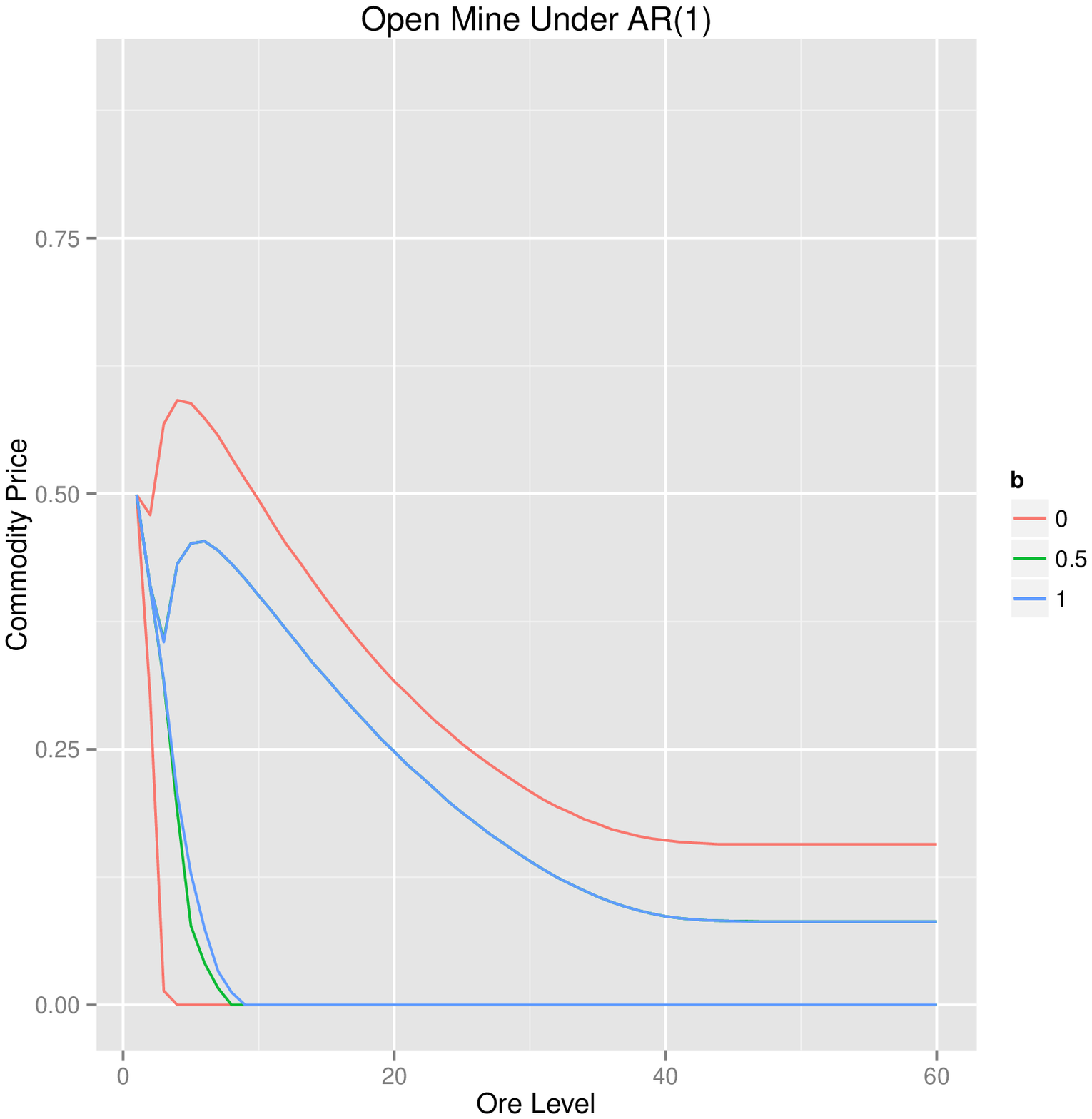}  
	\caption{Optimal policy at $t = 0$ under different values of $b = 0, 0.5, 1$.
		The policies under $\phi=1$ are depicted on the left while the right
		shows $\phi=0.6$. Here the blue and green curves are almost identical. \label{plot_contract}}
\end{figure}

The transition rule for $(Z_t)_{t=0}^T$ and the cashflow $h_t$ from
the previous part will be reused, as well as the grid and the distribution
sampling.  The results in Table \ref{table_contract} were generated
using $b = 0, 1$, $\phi = 1, 0.6$ and $w = 0$. We set
\[ p_t^* = 
\begin{cases} 
p_0^{(1)} - \frac{3}{4}(t-1) 
& \mbox{for $t = 5, 9, 13, \dots, 41 $};\\
p_0^{(1)} & \mbox{otherwise} \end{cases}\]
where $p_0^{(1)}$ is the starting commodity level. Under the above
specification of $p_t^*$, the resource is contracted to deliver 3 units of
the commodity by the first year, 6 units by the second year, ..., and
30 units by the end of the first decade. Using $K = 500$
sample paths and $I = 500$ sub-simulations in the diagnostics
 returns tight bounds and
small standard errors, which are given in the parenthesis.

In Table \ref{table_contract}, physical delivery obligations seem to
have little effect under strong mean reverting prices.  In contrast,
there is a significant impact under geometric Brownian motion
prices. This can also be seen in the optimal policies for an opened
resource at $t = 0$, as depicted in Figure \ref{plot_ar}.
\\

\begin{figure}[h!]
  \includegraphics[height=2.3in,width=0.5\textwidth]{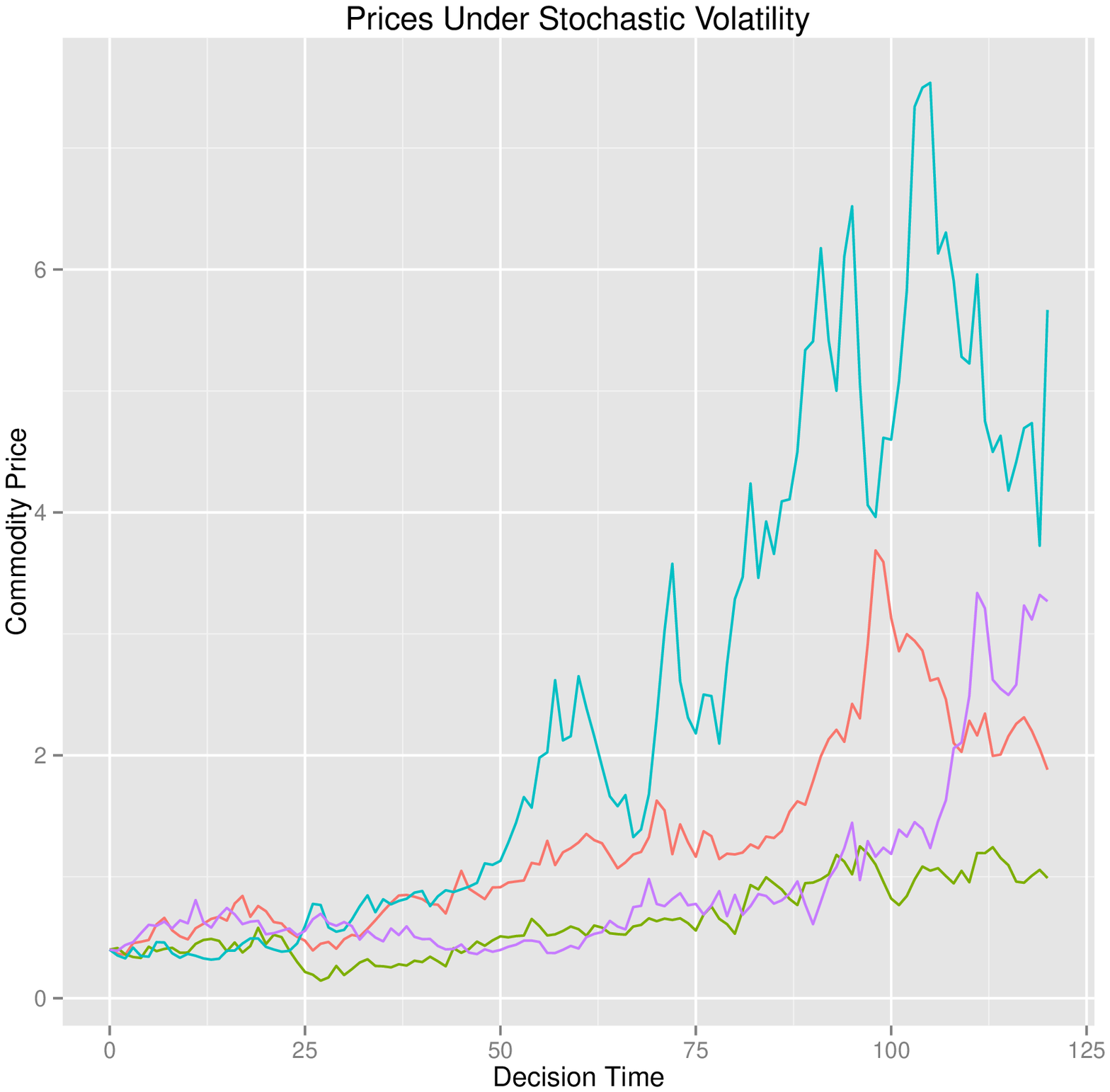}  
  \includegraphics[height=2.3in,width=0.5\textwidth]{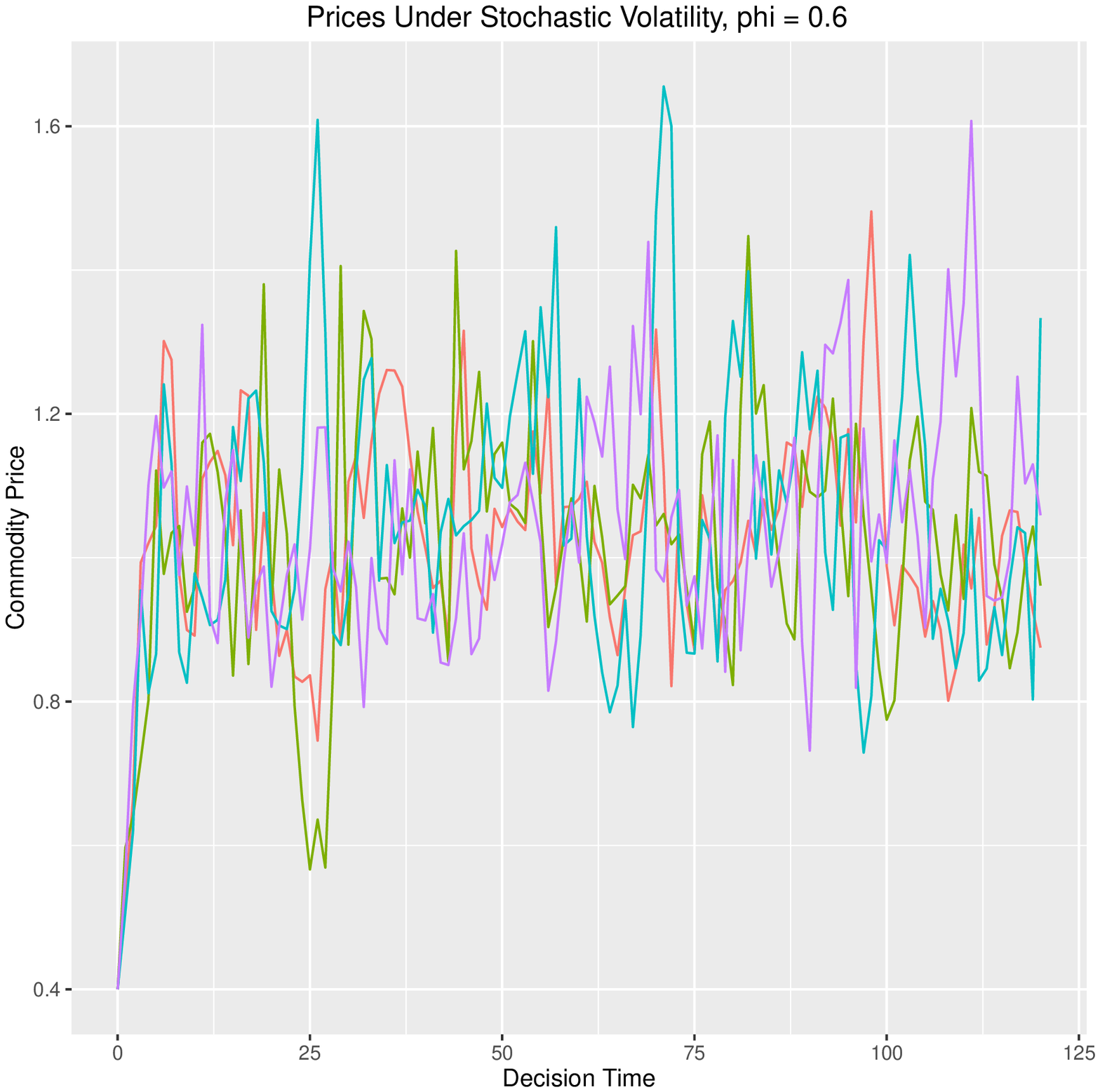}  
  \caption{Sample paths for $(\tilde Z_t)_{t=0}^T$ under no mean reversion (left)
    and mean reversion (right)\label{garchpaths}}
\end{figure}

\noindent
{\bf GARCH-type Process:} Such process are popular in modeling the
effects of time-changing fluctuation intensity, the so-called
volatility.  However, the original GARCH definition includes
non-linear recursion and is not covered by our approach which requires
a linear state dynamics.  For this reason, we suggest a simple
modification in order to retain this characteristic feature. In what
follows, we consider a GARCH(1,1)-like model. Having assumed the
controlled operation mode dynamics $(P_t)_{t=0}^T$ as in the original
geometric Brownian motion case and the absence of any delivery
contracts, we define the process $(\log \tilde Z_t)_{t=0}^T$
recursively using independent standard normally distributed random
variables $(N_t)_{t=1}^T$ as
\begin{align*}
  & \sigma^2_{t+1} = \sigma^2{\beta} 
    + \beta_1 \sigma^2_t + \beta_2 Y^2_t \\
  & Y_{t+1}^2 = \sigma^2_{t+1} N_t^2 \\
  & \log (\tilde Z_{t+1}) = \kappa\Delta + \phi \log (\tilde Z_{t}) +
    \sigma^2_{t+1} \sqrt{\Delta} N_{t+1}
\end{align*}
for $t=0, \dots, T-1$, with initial values $\sigma_0^2 =\sigma^2 =
\sqrt{0.08}$ $Y_0^2,\tilde Z_{0} \in \mathbb{R}_{+}$ and parameters
$\kappa = \mu - \sigma^4/2 = 0.05$, $\beta_1, \beta_2 \in \rr_+ $ with
$\beta_1+\beta_2 \in [0, 1]$ such that $\beta=1-\beta_1-\beta_2$. With
this definition, $(\sigma_t^2)_{t=0}^T$ follows the same recursion as
the conditional variance the original GARCH(1,1) process and can be
interpreted as the volatility proxy of the commodity price $(\tilde
Z_t)_{t=0}^T$ since $\sigma^2_{t+1} \sqrt{\Delta}$ is the conditional standard
deviation of the increment $\log \tilde Z_{t+1} - \log \tilde Z_t$ for
$t=0, \dots T-1$.  Positive values of $\beta_1, \beta_2$ induces
volatility clustering and a mean reversion is induced when $\phi \in
]0, 1[$.  With this choice, we define the linear state evolution $(
Z_t)_{t=0}^T$ for the state variables
$$ 
Z_{t} = [ \begin{array}{cccc} 1, & \sigma_{t}^2, & Y_{t}^2, & \log \tilde Z_{t}, \end{array} ]^\top, \qquad t=0, \dots T
$$
 in the required form $Z_{t+1}=W_{t+1}Z_t$ as
\begin{equation} \nonumber
Z_{t+1} = \begin{bmatrix} 1 & 0 & 0 & 0 \\ 
    \sigma^2 \beta & \beta_1 & \beta_2 & 0 \\
    \sigma^2 {\beta} N_{t+1}^2 & \beta_1 N_{t+1}^2 & \beta_2 N_{t+1}^2 & 0 \\
    \kappa\Delta + \sigma^2 \beta  \sqrt{\Delta} N_{t+1} &
    \beta_1 \sqrt{\Delta} N_{t+1} &  \beta_2 \sqrt{\Delta} N_{t+1} & \phi \end{bmatrix}
  \begin{bmatrix} 1 \\ \sigma_t^2 \\ Y_t^2 \\ \log \tilde Z_t \end{bmatrix}
\end{equation} 
Figure \ref{garchpaths} depicts the sample paths for the commodity
price $(\tilde Z_t)_{t=0}^T$ for parameters: $\tilde Z_0 = 0.4 $, $\phi = 1$
and $0.6$, $\kappa = 0.05$, $\beta_1 = 0.8$, $\beta_2 = 0.1$,
$\sigma^2 = \sigma_0^2 = \sqrt{0.08}$ and $Y_0^2 = 1$. The left
plot shows the sample paths when there is no mean reversion i.e.
$\phi = 1$ while the right contains mean reversion i.e.
$\phi = 0.6$. The behaviour of the prices are clearly very different.
Since the logarithmic commodity price is now contained in the fourth
component $z^{(4)}$ of the state vector $z=(z^{(i)})_{i=1}^4 $, we
adjust the cash flow function as
\[h_{t}(z) = 5\Delta \exp\left(z^{(4)}\right)\exp(-(r+\zeta)t\Delta) 
-2.5\Delta \exp ((\rho - r - \zeta)t\Delta).\]

In the previous examples we have used equally spaced grids. However,
for this four dimensional state space we generate the so-called
stochastic grid using clustering of a point cloud which we obtain from
several Monte Carlo runs of the sample paths simulations.  Such grid
generation procedure is appropriate for problems involving high
dimensional state spaces. The primal and dual values in Table
\ref{table_garch} were generated using a $2,000$ point stochastic grid
starting from $(1,\sqrt{0.08},1,\log (0.4))$, while the stochastic
grid itself was generated using $1,000$ sample paths. The primal and
dual values were generated using $K=500$ paths and $I=500$
sub-simulations. The Table \ref{table_garch} contains the upper and
lower bounds for the value of the commodity reserve for $\beta_1 =
0.8$, $\beta_2 = 0.1$ and $Y_0^2 = 1$.  These values were computed
using $10,000$ point disturbances discretization based on equidistant
quantiles of the standard normal.

\begin{table}[h!]
  \centering
  \setlength\tabcolsep{6pt}
  \caption{Resource Valuation Under Stochastic Volatility \label{table_garch}}
  \begin{tabular}{lrrrrr}
    & &\multicolumn{2}{c}{Opened Resource}&\multicolumn{2}{c}{Closed Resource}\\
    $\phi$ & $\tilde Z_0$ & Primal & Dual & Primal & Dual \\
    \hline
      1 & 0.3 &   3.4337(.1069) &  4.4240(.4663)&  3.6337(.1069) &  4.6240(.4663) \\
        & 0.4 &   7.5215(.2035) &  8.8880(.6060)&  7.7215(.2035) &  9.0878(.6060) \\
        & 0.5 &  12.7777(.6901) & 13.8656(.7422)& 12.9777(.6901) & 14.0614(.7422) \\
    \hline
    0.8 & 0.3 &   6.5229(.0109) &  6.5372(.0152)&  6.5232(.0110) &  6.5376(.0152) \\
        & 0.4 &   7.4872(.0112) &  7.5005(.0153)&  7.3028(.0112) &  7.3166(.0153) \\
        & 0.5 &   8.3380(.0120) &  8.3512(.0159)&  8.1380(.0120) &  8.1513(.0159) \\
    \hline
    0.6 & 0.3 &   7.7928(.0052) &  7.7960(.0058)&  7.7262(.0052) &  7.7294(.0059) \\
        & 0.4 &   8.2726(.0051) &  8.2759(.0057)&  8.0755(.0051) &  8.0788(.0057) \\
        & 0.5 &   8.6992(.0051) &  8.7024(.0058)&  8.4992(.0051) &  8.5024(.0058) \\
  \end{tabular}
  \caption*{Resource valuation where $\kappa = 0.05$, $\beta_1 = 0.8$, $\beta_2 = 0.1$, 
    $\sigma_0^2 = \sqrt{0.08}$ and $Y_0^2 = 1$. The standard errors are given in the paranthesis.}
\end{table}

Table \ref{table_garch} also clearly shows that bounds are less
accurate for higher values of $\phi$ due to lack of concentration due
to weaker mean reversion.  A more dense grid and a larger sample of
the disturbances should yield more accurate results.

\section{Conclusion \label{sec-end}}
This paper represents a first application of pathwise approach to
dynamic programming in the area of commodity resource valuation and
extraction.  Our study shows that high-quality solutions can be
obtained with low computational efforts for realistic applications.
By representing the extraction problem as a convex stochastic
switching problem in discrete time, this paper has demonstrated that
accurate and helpful practical insights can be deduced from
approximate numerical solutions.  Due to a standardized problem
formulations, other situations can easily be adopted and our concepts
can generalized and extended to other application areas of optimal
stochastic switching.

\bibliography{resource}
\bibliographystyle{plain}
\end{document}